\documentclass{article}
\usepackage[T1]{fontenc}
\usepackage{lmodern}
\input%
{glyphtounicode}
\usepackage[letterpaper,margin=1in]{geometry}
\usepackage{mathtools, amsfonts, amsthm, amssymb, graphicx, mleftright}
\usepackage{enumitem}
\usepackage[colorlinks,hypertexnames=false]{hyperref}
\usepackage{url}

\numberwithin{equation}{section}

\newtheorem{thm}{Theorem}[section]
\newtheorem{lem}[thm]{Lemma}
\newtheorem{exer}[thm]{Exercise}
\newtheorem{prop}[thm]{Proposition}
\newtheorem{conj}[thm]{Conjecture}
\theoremstyle{definition}
\newtheorem{qn}[thm]{Question}
\theoremstyle{remark}
\newtheorem{rem}[thm]{Remark}

\DeclareMathOperator{\lcm}{lcm}
\DeclareMathOperator{\sgn}{sgn}
\DeclareMathOperator{\Mat}{Mat}

\newcommand{\GL}{\mathrm{GL}}
\newcommand{\PGL}{\mathrm{PGL}}

\newcommand{\QQ}{\mathbb{Q}}
\newcommand{\PP}{\mathbb{P}}
\newcommand{\RR}{\mathbb{R}}
\newcommand{\ZZ}{\mathbb{Z}}

\newcommand{\V}{\mathcal{V}}
\newcommand{\W}{\mathcal{W}}

\newcommand{\longto}{\mathop{\longrightarrow}\limits}
\newcommand{\textand}{\quad \text{and} \quad}
\renewcommand{\to}{\mathop{\rightarrow}\limits}

\newcommand{\floor}[1]{\left\lfloor #1 \right\rfloor}
\newcommand{\ceil}[1]{\left\lceil #1 \right\rceil}

\mleftright
\renewcommand{\epsilon}{\varepsilon}

\title{Two-point Pad\'e approximation\\
and patterns in continued fractions for $(F_{n+1}/F_n)^k$
}

\author{Evan M. O'Dorney}

\begin{document}

\maketitle

\begin{abstract}
  We explore the simple continued fraction expansion for the $k$th power of the Fibonacci convergent $F_{n+1}/F_n$, as $n$ grows. For most $k$, this continued fraction consists of $k$ linearly growing \emph{patches} of infinite periodic continued fractions separated by \emph{junctions} that vary periodically in $n$. But under special conditions, two patches merge or, more strikingly, an exponentially growing \emph{monster term} replaces two patches and the junction between them.

  We trace these phenomena to properties of the Pad\'e approximants that approximate $x^k$ near $x = \phi$ and $x = \bar\phi$ simultaneously. We prove the patterns for second-order unimodular recurrences and compute several more instances of this rare and apparently previously unnoticed monster term.
\end{abstract}

\section{Introduction}
In 2020, Chu \cite{Chu2020} asked which of the pair of linear Diophantine equations
\begin{equation} \label{eq:frob}
  a x + b y = \frac{a b - a - b \pm 1}{2}, \qquad x \geq 0, \quad y \geq 0
\end{equation}
has an integer solution $(x,y)$ for given positive coprime coefficients $a$ and $b$, a question arising out of studying the middle coefficient of cyclotomic polynomials. It turns out that \eqref{eq:frob} is solvable for exactly one of the two signs (this is an instance of Frobenius symmetry of two-generator numerical semigroups); and in the case that $(a,b) = (F_n, F_{n+1})$ are consecutive Fibonacci numbers, Chu observed and proved that the sign, and associated formulas for $x$ and $y$, oscillate with period $6$ in $n$. This discovery inspired a succession of Polymath Jr.\ projects \cite{ChenEtAl2025,ArachchiEtAl2026,ChuGulecha2026} in which Chu, Miller, and their mostly undergraduate collaborators analyzed the solutions to \eqref{eq:frob} for $a, b$ drawn from an assortment of number sequences, including squares and cubes of Fibonacci numbers, where similar periodicities appear.

In this paper, rather than studying \eqref{eq:frob} directly, we step back to the object that computes the solution to $a x + b y = c$ in general: the continued fraction expansion of $b/a = F_{n+1}^k/F_n^k$. It is well known that the Fibonacci ratios
\begin{equation} \label{eq:fib_1}
  \frac{F_{n+1}}{F_n} = \overbrace{[1,1,\ldots,1]}^n
\end{equation}
are the convergents to the periodic continued fraction $[\overline{1}]$ for the golden ratio $\phi = (1 + \sqrt{5})/2$. Raising them to small powers yields repeating patterns:
\begin{equation} \label{eq:fib_2}
  \begin{aligned}
    (F_2/F_1)^2 &= [1] \\
    (F_3/F_2)^2 &= [4] \\
    (F_4/F_3)^2 &= [\mathbf{2}, \mathbf{3}, 1] \\
    (F_5/F_4)^2 &= [\mathbf{2}, 1, \mathbf{3}, 1, 1] \\
    (F_6/F_5)^2 &= [\mathbf{2}, 1, 1, \mathbf{3}, 1, 1, 1] \\
    (F_7/F_6)^2 &= [\mathbf{2}, 1, 1, 1, \mathbf{3}, 1, 1, 1, 1] \\
    (F_8/F_7)^2 &= [\mathbf{2}, 1, 1, 1, 1, \mathbf{3}, 1, 1, 1, 1, 1].
  \end{aligned}
\end{equation}
Excluding the first two values (which can be shoehorned in via a suitable interpretation of negative repetitions of a continued fraction term), the pattern continues with one more $1$ on each side of the $3$. The pattern for cubes is only slightly more complex:
\begin{equation}\label{eq:fib_3}
  \begin{aligned}
    (F_9/F_8)^3 &= [4,4,\mathbf{10},4,\mathbf{2},\mathbf{2},1,1,1,1,1] \\
    (F_{10}/F_9)^3 &= [4,4,\mathbf{3},\mathbf{2},\mathbf{3},4,\mathbf{2},\mathbf{2},1,1,1,1,1,1] \\
    (F_{11}/F_{10})^3 &= [4,4,4,\mathbf{1},\mathbf{1},\mathbf{1},4,4,\mathbf{2},\mathbf{2},1,1,1,1,1,1,1] \\
    (F_{12}/F_{11})^3 &= [4,4,4,\mathbf{10},4,4,\mathbf{2},\mathbf{2},1,1,1,1,1,1,1,1] \\
    (F_{13}/F_{12})^3 &= [4,4,4,\mathbf{3},\mathbf{2},\mathbf{3},4,4,\mathbf{2},\mathbf{2},1,1,1,1,1,1,1,1,1] \\
    (F_{14}/F_{13})^3 &= [4,4,4,4,\mathbf{1},\mathbf{1},\mathbf{1},4,4,4,\mathbf{2},\mathbf{2},1,1,1,1,1,1,1,1,1,1].
  \end{aligned}
\end{equation}
Observe that the continued fraction for each $(F_{n+1}/F_n)^3$ has the structure
\begin{equation}
  \mathbf{P}_0^{\floor{(n-1)/3}} \mathbf{J}_1 \mathbf{P}_1^{\floor{(n-4)/3}} \mathbf{J}_2 \mathbf{P}_2^{n-3}
\end{equation}
consisting of three \emph{patches} $\mathbf{P}_i^{e_i(n)}$ of repeated $4$'s and $1$'s separated by \emph{junctions} $\mathbf{J}_i$: Junction $\mathbf{J}_1$ is alternately $[10]$, $[3, 2, 3]$, and $[1,1,1]$, while $\mathbf{J}_2 = [2, 2]$.

As the exponent $k$ increases, we start seeing more and more junction types as well as more complex patches: fourth powers give
\begin{equation}\label{eq:fib_4}
  \begin{aligned}
    (F_{26}/F_{25})^4 &= [\mathbf{6},1,5,1,5,1,5,1,5,1,5,1,\mathbf{1},3,1,2,3,2,1,3,33,\mathbf{4},\mathbf{6},11^9] \\
    (F_{27}/F_{26})^4 &= [\mathbf{6},1,5,1,5,1,5,1,5,1,5,1,\mathbf{37},3,1,2,3,2,1,3,\mathbf{31},\mathbf{1},\mathbf{9},11^9] \\
    (F_{28}/F_{27})^4 &= [\mathbf{6},1,5,1,5,1,5,1,5,1,5,1,\mathbf{4},\mathbf{4},33,3,1,2,3,2,1,3,\mathbf{10},\mathbf{1},\mathbf{10},11^9] \\
    (F_{29}/F_{28})^4 &= [\mathbf{6},1,5,1,5,1,5,1,5,1,5,1,\mathbf{6},\mathbf{1},3,2,1,3,33,3,1,2,3,2,1,\mathbf{15},11^{10}] \\
    (F_{30}/F_{29})^4 &= [\mathbf{6},1,5,1,5,1,5,1,5,1,5,1,5,1,\mathbf{1},3,1,2,3,2,1,3,33,3,1,2,\mathbf{2},\mathbf{1},\mathbf{1},11^{11}] \\
    (F_{31}/F_{30})^4 &= [\mathbf{6},1,5,1,5,1,5,1,5,1,5,1,5,1,\mathbf{37},3,1,2,3,2,1,3,33,\mathbf{4},\mathbf{6},11^{11}] \\
    (F_{32}/F_{31})^4 &= [\mathbf{6},1,5,1,5,1,5,1,5,1,5,1,5,1,\mathbf{4},\mathbf{4},33,3,1,2,3,2,1,3,\mathbf{31},\mathbf{1},\mathbf{9},11^{11}].
  \end{aligned}
\end{equation}
(Within the continued-fraction brackets, $a^b$ means that the term $a$ is repeated $b$ times.)  Here one patch consists of segments of the cyclic word $[3, 1, 2, 3, 2, 1, 3, 33]$, and if we extend the table further, we will see that this patch widens by one full period for every $20$ values of $n$. Not coincidentally, this word is the periodic portion of the continued fraction
\begin{equation}
  15 \phi = [24, \overline{3, 1, 2, 3, 2, 1, 3, 33}]
\end{equation}
of an integer multiple of the golden ratio.

The values $k = 6, 7, 8, \ldots$ show similar behavior, with increasing complexity. However, by far the most striking behavior is that for $k = 5$:
\begin{equation}\label{eq:fib_5_table}
  \begin{aligned}
    (F_{21}/F_{20})^5 &= [11^3,10,1,1,\mathbf{33116027},1,17,11^2,10,1,4,11^3] \\
    (F_{22}/F_{21})^5 &= [11^4,2,1,\mathbf{86698886},2,5,11^3,21,11^3] \\
    (F_{23}/F_{22})^5 &= [11^4,\mathbf{5674515856},1^{19},2,9,11^3] \\
    (F_{24}/F_{23})^5 &= [11^4,7,1,\mathbf{594243013},2,5,11^4,1,11^4] \\
    (F_{25}/F_{24})^5 &= [11^4,12,1,\mathbf{1555748409},6,1,10,11^3,4,1,10,11^3] \\
    (F_{26}/F_{25})^5 &= [11^4,10,2,\mathbf{4073002214},18,11^3,10,1,4,11^4] \\
    (F_{27}/F_{26})^5 &= [11^5,3,\mathbf{10663258234},1,1,5,11^4,21,11^4] \\
    (F_{28}/F_{27})^5 &= [11^4,10,1,\mathbf{697919312217},2,1^{22},2,9,11^4] \\
    (F_{29}/F_{28})^5 &= [11^5,8,\mathbf{73087059232},1,1,5,11^5,1,11^5] \\
    (F_{30}/F_{29})^5 &= [11^5,13,\mathbf{191344405207},1,5,1,10,11^4,4,1,10,11^4].
  \end{aligned}
\end{equation}
The arresting feature is the appearance of one
``monster'' term in the continued fraction, a term of order $\phi^{2n}$ (with a leading factor depending on $n$ modulo $5$).

A related phenomenon was noted by French in 2006 \cite{French2006}:
the continued fraction of $\sqrt[5]{F_{n+5}/F_n}$ has a large partial
quotient $F_{2n+5} - 1 - 3(-1)^n$. Underlying both large terms is the close rational approximation
\begin{equation} \label{eq:fib_5_approx}
  \left(\frac{F_{n+1}}{F_n}\right)^5 = \frac{F_{n+3}}{F_{n-2}} + O(F_n^{-4}).
\end{equation}
Here, an error of approximation less than $1/(2 F_{n-2}^2)$ is enough to guarantee that the fraction on the right-hand side is a convergent; further accuracy implies a correspondingly large term.

In this paper, we ask what happens when we vary the exponent $k$ as well as the Fibonacci-type recurrence under consideration.
\begin{qn}\label{qn:main}
  Let $\{u_n\}$ be a sequence of nonnegative integers given by a recurrence relation
  \begin{equation}\label{eq:u_n}
    u_{n+2} = a u_{n+1} + b u_n, \qquad u_0 = 0, \quad u_1 = 1,
  \end{equation}
  and let $k$ be a positive integer.
  \begin{enumerate}
    \item\label{qn:JPJ} When does the continued fraction for $(u_{n+1}/u_n)^k$ display a structure of patches and junctions for $n$ large enough?
    \item\label{qn:bound} Can we bound the periods involved?
    \item\label{qn:monster} When is there a monster term?
  \end{enumerate}
\end{qn}

It is too much to expect orderly results for all recurrence parameters $a$ and $b$. For instance, taking $a = 2$, $b = 5$, we have, already for $k = 1$:
\begin{equation} \label{eq:chaos}
  \begin{aligned}
    u_{10}/u_{9}  &= [3; \mathbf{2, 4}, 3, 1, 4, 1, 20, 1, 2]                    \\
    u_{11}/u_{10} &= [3; \mathbf{2, 4}, 1, 1, 7, 107, 3]                         \\
    u_{12}/u_{11} &= [3; \mathbf{2, 4, 2}, 2, 2, 4, 21, 3, 1, 3]                 \\
    u_{13}/u_{12} &= [3; \mathbf{2, 4, 2}, 6, 2, 5, 1, 3, 2, 4, 1, 6]            \\
    u_{14}/u_{13} &= [3; \mathbf{2, 4, 2}, 3, 1, 10, 1, 6, 2, 4, 3, 2, 1, 2]     \\
    u_{15}/u_{14} &= [3; \mathbf{2, 4, 2, 4}, 1, 2, 1, 1, 36, 4, 5, 1, 1, 1, 3]  \\
    u_{16}/u_{15} &= [3; \mathbf{2, 4, 2, 4, 2}, 1, 8, 7, 4, 29, 5, 2]           \\
    u_{17}/u_{16} &= [3; \mathbf{2, 4, 2, 4, 2}, 43, 1, 2, 4, 2, 5, 2, 3, 2, 3].
  \end{aligned}
\end{equation}
As shown in bold, there is an initial periodic patch converging to the limiting value $\alpha = [3; \overline{2,4}]$, but the rest of the terms are apparently random. Of course, \emph{proving} that the terms obey the Gauss--Kuzmin distribution (which governs the continued fraction of a random real number) is not necessarily easy.

We therefore restrict to the case $b = \pm 1$ for the bulk of the paper. However, the necessary and sufficient condition for most (but not all) of the results is $b \mid a^2$, as we will discuss at the end of the paper.

\begin{thm}\label{thm:main_intro}
  There are polynomials $\Delta_{k,i}(a,b) \in \ZZ[a,b]$ for integers $0 \leq i \leq k$, with $\Delta_{k,0} = \Delta_{k,k} = 1$, with the following properties.

  Consider a recurrence
  \begin{equation*}
    u_0 = 0, \quad u_1 = 1, \quad u_{n+2} = a u_{n+1} + b u_n
  \end{equation*}
  where $b = \pm 1$ and $a \geq 2 - b$, so that the characteristic polynomial $h(x) = x^2 - a x - b$ has irrational real roots $\alpha > \bar\alpha$. Then for each $k \geq 1$,
  \begin{enumerate}[label=\textup{(\alph*)}]
    \item\label{main:monster} If $\Delta_{k,i}(a,b) = 0$ for some $i$, then for all sufficiently large $n$, the power $(u_{n+1}/u_n)^k$ has a continued fraction with a ``monster term\kern-1pt'' $M$ bounded below and above by constant multiples of $\alpha^{2n}$.
    \item\label{main:JPJ} If, on the other hand, $\Delta_{k,i}(a,b) \neq 0$ for all $i$, let $m_i$ be the least positive integer for which
    \begin{equation*}
      \Delta_{k,i}(a,b) \mid u_{m_i},
    \end{equation*}
    and let $g_i = \gcd(m_i, m_{i+1})$.

    Then for $n$ sufficiently large, $(u_{n+1}/u_n)^k$ has a continued fraction with a patches-and-junctions decomposition
    \begin{equation}
      \mathbf{J}_{0,0} \mathbf{P}_{0,0}^{e_0(n)} \mathbf{J}_{1, n \bmod m_1} \mathbf{P}_{1, n \bmod g_1}^{e_1(n)} \mathbf{J}_{2, n \bmod m_2} \cdots \mathbf{J}_{k-1, n \bmod m_{k-1}} \mathbf{P}_{k-1, 0}^{e_{k-1}(n)}
    \end{equation}
    where the continued-fraction words $\mathbf{P}_{i,n}$, $\mathbf{J}_{i,n}$ depend only on $n$ modulo $g_i$ and $m_i$ respectively, and the number $e_i(n)$ of repetitions within each patch is a linear function of $n$ when $n$ is in a fixed congruence class modulo $\lcm\{m_i, m_{i+1}\}$. In particular, the terms of the continued fraction remain bounded as $n$ grows.
  \end{enumerate}
\end{thm}

\subsection{Why study Fibonacci numbers in the age of AI?}

We live in a moment of great uncertainty about the status of mathematics and mathematicians. This article is intended for a broad audience, including students (and those who supervise students) who are drawn by the beauty of mathematics but are doubting whether it is worthwhile to persevere in learning ever more involved notions and theories. Here is a word to them.

The history of mathematics has seen the rise and fall of many specialized fields of study. Euclidean geometry---the very language of ancient Greek mathematics---has been reduced to a training ground in elementary logic and a catalog of theorems rarely cited beyond high-school olympiads. Nineteenth-century mathematicians filled massive tomes with properties of special functions (elliptic functions, orthogonal polynomials, etc.). By the mid-20th century this field was moribund, increasingly regarded as a ``zoo'' of identities to be replaced by Bourbaki-style abstract theory \cite{Beals2012PDE};  ironically, special functions experienced a resurgence in the early 21st century via the rise of symbolic computation \cite{Wolfram2005SpecialFunctions}.
The decline of a field is pervaded by a sense, not that every possible result therein has been discovered (indeed, in many fields, there is an easy infinite supply of theorems), but that anything new would be of necessity a recombination of old ideas and methods. We are now facing a threat of this sort, not just to individual fields but to the whole enterprise of mathematics: if a machine can do our research for us---pose problems, try proof techniques, write papers---does it not, in some sense, expose us all as zookeepers collecting specimens ever deeper in a fixed (though infinite) tree of possible research paths?

This paper studies the properties of Fibonacci numbers and allied recurrences, which, long before the founding of the \emph{Fibonacci Quarterly} in 1963, had acquired the air of a zoo (and a zoo far more recreational than the special functions that solve industrially significant differential equations!). The proofs in this paper are concatenations of techniques at least a century old. Given any proposed identity in Fibonacci numbers, a modern chatbot can produce a correct proof or counterexample within minutes. However, this is not a reason to stop investigating the Fibonacci sequence, or mathematics at large.

Mathematics was never solely about proving or disproving conjectures; and the purpose of this paper is not, primarily, to prove new theorems. Rather, I present it as a small gift in the \emph{formation of the mathematical sense,} even of the \emph{artistic sense} within mathematics. The questions ``What objects should we study? What operations should we apply to them? What results are noteworthy if true?'' are fundamentally not reducible to mechanical calculation. The danger of AI to mathematics is not so much that it will do our research for us; it is that too many people will stop caring whether research gets done at all. And that is not such a happy future. Already, it is becoming clear that the best workers in any field, assisted by AI, are those with track records as good workers before AI (often coupled with the adaptability to incorporate AI into the workflow). This includes valuable work in pure and applied mathematics. Let works like this one keep the flame of interest alive!

\subsection{Conventions}
In this paper, all continued fractions will be simple continued fractions $[a_0; a_1, \ldots, a_n]$ or $[a_0; a_1, a_2, \ldots]$ with positive integer terms $a_i$, $i \geq 1$. We allow $a_0 \in \ZZ$. We do \emph{not} assume that the last term $a_n$ is at least $2$: each rational number has two continued fraction expansions, and we reserve the right to display the choice yielding a more elegant pattern, as we have been doing since the first page of the article.

This paper is peppered with \emph{Remarks} and \emph{Exercises} which contain statements whose proofs are omitted. The statements are meant to give a fuller picture of the phenomena in question; the proofs are not meant to be trivial, but with an AI tool the reader will not have any trouble reconstructing them. The results in this paper are not yet formalized; formalization would be itself a reasonable exercise.

\subsection{AI disclosure}
All mathematical proofs in this article are based on work generated by ChatGPT 6, entirely rewritten by the author.

\section{Pad\'e approximations and monster terms} \label{sec:pade_monster}
One of the key tools for explaining the patterns observed above is \emph{Pad\'e approximation}, a version of Diophantine approximation for functions which we briefly recall here. Classically, Pad\'e approximation seeks to express
a given power series $f(x)=a_0+a_1x+a_2x^2+\cdots$ to as high order as
possible by a rational function $P(x)/Q(x)$ with restricted degrees
$\deg P\leq m$, $\deg Q\leq n$. For generic $f$, one can achieve contact order
\begin{equation*}
  f(x) = \frac{P(x)}{Q(x)} + O(x^{m+n+1});
\end{equation*}
the approximant $P/Q$ is unique if it exists
\cite{WikipediaPade,Baker2012}. The Pad\'e approximants are the convergents of a certain nonsimple continued fraction, and in many cases they appear in the corresponding simple continued fraction. For instance, in classical expansions such as
\begin{equation}
  e^{\frac{1}{1000}} = [1; \mathbf{999}, 1, 1, \mathbf{2999}, 1, 1, \mathbf{4999}, 1, 1, \mathbf{6999}, \ldots]
\end{equation}
(see \cite{Cohn06,McCabe09}), as well as more fleeting patterns like
\begin{equation}
  \sqrt{10000^2 + 3} = [10000; \mathbf{6666}, 1, 2, \mathbf{2221}, 1, 8, \mathbf{740}, 1, 1, 1, 2, 2, 1, \mathbf{246}, 4, 1, 3, 4, \mathbf{82}, 14, 3, 2, 2, \mathbf{27}, \ldots],
\end{equation}
the large terms shown in bold measure the accuracy of successive Pad\'e approximants to the appropriate function.

The continued fractions that we study here also relate to Pad\'e approximations, but of a less common \emph{two-point} type. For instance,  \begin{equation} \label{eq:sol_sq}
  \left(\frac{F_{n+1}}{F_n}\right)^2 = \frac{F_{n}}{F_{n-2}} + O(F_n^{-2}) = \cfrac{1}{2 - \cfrac{F_{n+1}}{F_n}} + O(F_n^{-2})
\end{equation}
has the shape $x^k \approx \frac{P(x)}{Q(x)}$, where $x = \frac{F_{n+1}}{F_n}$; but here the approximation is not centered at $x = 0$, but at $x = \phi$ and (considering $n \to -\infty$) also at the algebraic conjugate $x = \bar\phi = (1 - \sqrt{5})/2$. This leads us to \emph{two-point Pad\'e approximation} \cite{Sidi1980}
where the aim is to find $P(x), Q(x)$ of restricted degrees for which
\begin{equation} \label{eq:Pade1}
  x^k \equiv \frac{P(x)}{Q(x)} \mod (x - \alpha)^i (x - \bar\alpha)^j.
\end{equation}
In our setting, $\alpha, \bar\alpha$ are the two roots of the characteristic polynomial $h$. If $P$ and $Q$ are to have integer coefficients, the contact orders $i, j$ may be taken equal, and \eqref{eq:Pade1} can be written as
\begin{equation} \label{eq:PadeE}
  x^k Q(x) - P(x) = E(x) h(x)^i
\end{equation}
for some $E(x) \in \ZZ[x]$ which represents the \emph{error} of approximation. Suppose $\deg P = \deg Q = d$. Since $P$ and $Q$ together have $2d + 2$ coefficients and the divisibility \eqref{eq:PadeE} imposes $2i$ homogeneous linear conditions, we \emph{usually} must go up to degree $d = i$ to get a nonzero solution $(P,Q)$; but occasionally an approximant of degree less than $i$ is \emph{unusually accurate:}
\begin{equation}\label{eq:fib_5_error}
  x^5 = \frac{2x+1 - (x^2+1)(x^2-x-1)^2}{2 - x} \equiv \frac{2x + 1}{2 - x} \mod (x^2 - x - 1)^2
\end{equation}
is a linear approximation achieving quadratic accuracy. Substituting $x = F_{n+1}/F_n$, this is none other than the approximation \eqref{eq:fib_5_approx} underlying the Fibonacci monster term for $k = 5$. (Classically, this phenomenon is called a ``Pad\'e defect,'' but we are reluctant to apply the pejorative label ``defect'' to a rare coincidence with exciting consequences!)

An unusually accurate Pad\'e approximant exists precisely when the linear system packaged in \eqref{eq:PadeE} is singular. In the present instance, the relevant matrix can be cast in a pleasantly simple form:
\begin{lem} \label{lem:Pade_Toeplitz}
  Let $0 \leq i \leq k$, and write
  \begin{equation*}
    h(x)^i = c_{0}x^{2i} + c_1 x^{2i-1} + \cdots + c_{2i-1} x + c_{2i}.
  \end{equation*}
  Form the $(k-i) \times (k-i)$ Toeplitz matrix
  \begin{equation}
    \mathbf{C}_{k,i} = \begin{bmatrix}
      c_i & c_{i+1} & \cdots & c_{k-1} \\
      c_{i-1} & c_i & \cdots & c_{k-2} \\
      \vdots & \vdots & \ddots & \vdots \\
      c_{2i-k+1} & c_{2i-k+2} & \cdots & c_i
    \end{bmatrix}
  \end{equation}
  of central coefficients of $h(x)^i$. (Any $c_j$ occurring with $j < 0$ or $j > 2i$ are interpreted as $0$.)
  The solutions to \eqref{eq:PadeE} of unusually low degrees
  \begin{equation*}
    \deg P < i, \qquad \deg Q < i, \qquad \deg E < k - i
  \end{equation*}
  are given by $E(x) = e_0 x^{k-i-1} + \cdots + e_{k-i-2} x + e_{k-i-1}$, where $\mathbf{e} = (e_{k-i-1}, \ldots, e_0)^\top$ satisfies the homogeneous linear system $\mathbf{C}_{k,i} \mathbf{e} = 0$. In particular, such a solution exists exactly when $\Delta_{k,i} = \det \mathbf{C}_{k,i}$ vanishes.
\end{lem}
\begin{proof}
  As $Q$ and $P$ vary, the left-hand side of \eqref{eq:PadeE} can be any integer combination of
  \begin{equation*}
    x^0, x^1, \ldots, x^{i-1}, \qquad x^{k}, x^{k+1}, \ldots, x^{k+i-1},
  \end{equation*}
  but the $k - i$ intermediate monomials
  \begin{equation*}
    x^i, x^{i+1}, \ldots, x^{k-1}
  \end{equation*}
  must appear with coefficient $0$. Computing the coefficients of these monomials in $E(x) h(x)^i$ yields the claimed linear system.
\end{proof}
For fixed $k$ and $i$, $\Delta_{k,i} = \det \mathbf{C}_{k,i}$ is a polynomial in $a$ and $b$. Note that $\Delta_{k,0} = \Delta_{k,k} = 1$. For a telling example,
\begin{equation} \label{eq:Delta_factored}
  \begin{gathered}
    \Delta_{5,2}(a,b) = \begin{vmatrix}
      a^2 - 2b & 2ab & b^2 \\
      -2a & a^2 - 2b & 2ab \\
      1 & -2a & a^2 - 2b
    \end{vmatrix} = a^6 + 2 a^4 b + 3 a^2 b^2 - 6 b^3
    = (a^2-b)(a^4+3a^2b+6b^2),
  \end{gathered}
\end{equation}
and the vanishing of the factor $a^2 - b$ at $a = b = 1$ corresponds to the Pad\'e approximant \eqref{eq:fib_5_error} and thus to the monster terms in \eqref{eq:fib_5_table}. The homogeneity is already noticeable:
\begin{lem}\label{lem:Delta_homog}
  The Pad\'e--Toeplitz determinant $\Delta_{k,i}$ has the weighted homogeneous form
  \begin{equation}\label{eq:Delta_cases}
    \Delta_{k,i}(a,b) = \begin{cases*}
      b^{\frac{i(k-i)}{2}} \delta_{k,i}(a^2/b), & $i(k-i)$ even \\
      - a b^{\frac{i(k-i) - 1}{2}} \delta_{k,i}(a^2/b), & $i(k-i)$ odd
    \end{cases*}
  \end{equation}
  where $\delta$ is a monic polynomial of degree $\floor{i(k-i)/2}$. Consequently, for $a, b \neq 0$, we can have $\Delta_{k,i} = 0$ only if $b \mid a^2$.
\end{lem}
\begin{proof}
  Give the polynomial ring $\ZZ[a,b]$ the grading where $\deg a = 1$ and $\deg b = 2$. Then each $c_j$ is homogeneous of degree $j$, so $\Delta_{k,i}$ is homogeneous of degree $i(k-i)$. This shows that no terms can occur other than the ones appearing on the right-hand side of \eqref{eq:Delta_cases}.

  The leading coefficient of $\delta$ is the coefficient in $\Delta_{k,i}$ of $a^{i(k-i)}$, the unique term with no factor of $b$. Setting $b = 0$, we observe that $h(x)^i = (x^2 - ax)^i$ has no terms of degree below $i$, so the Toeplitz matrix $\mathbf{C}_{k,i}$ is lower triangular and $\Delta_{k,i} = (-a)^{i(k-i)}$.

  Consequently, by the rational root test, $\Delta_{k,i}$ can vanish only if $b \mid a^2$.
\end{proof}

\subsection{No ``super-sized monsters''}

The next lemma shows that the coincidence cannot be \emph{too} strong.
\begin{lem}\label{lem:no_super}
  Two consecutive Pad\'e--Toeplitz determinants $\Delta_{k,i}, \Delta_{k,i+1}$ cannot simultaneously vanish for $a, b \neq 0$, and a vanishing $\Delta_{k,i} = 0$ corresponds to a Pad\'e approximant $P(x)/Q(x)$ with $P, Q$ coprime rational polynomials of degree exactly $i - 1$.
\end{lem}
\begin{proof}
  Suppose $\Delta_{k,i} = 0$ for some $0 < i < k$, and assume that $\Delta_{k,i-1} \neq 0$, that is, we are at the left end of a maximal block of consecutive vanishing Pad\'e--Toeplitz determinants. Then there exist rational polynomials $P,Q$ of degree at most $i - 1$ such that
  \begin{equation} \label{eq:super_1}
    P(x) \equiv x^k Q(x) \mod h(x)^i.
  \end{equation}
  If $P$ and $Q$ have a common monic irreducible factor in $\QQ[x]$, which may or may not be $h$, then dividing it out gives a solution to $\widetilde{P}(x) \equiv x^k \widetilde{Q}(x) \mod h(x)^{i-1}$ in polynomials of degree at most $i - 2$, implying $\Delta_{k,i-1} = 0$, a contradiction. So $P$ and $Q$ are coprime in $\QQ[x]$, and in particular, $Q(\alpha)$ and $P(\alpha) = \alpha^k Q(\alpha)$ are nonzero.

  Taking the logarithmic derivative of both sides of \eqref{eq:super_1}, which reduces the order of contact by at most $1$, we get
  \begin{alignat*}{2}
    \frac{P'(x)}{P(x)} &\equiv \frac{k}{x} + \frac{Q'(x)}{Q(x)} &&\mod h(x)^{i-1} \\
    x P'(x) Q(x) - x P(x) Q'(x) - k P(x) Q(x) &\equiv 0 &&\mod h(x)^{i-1}.
  \end{alignat*}
  The left-hand side is a polynomial of degree at most $2i - 2$, and $h(x)^{i-1}$ also has degree $2 i - 2$. So there is a constant $c$ for which
  \begin{equation} \label{eq:xPQc}
    x P'(x) Q(x) - x P(x) Q'(x) - k P(x) Q(x) = c h(x)^{i-1}.
  \end{equation}
  We cannot have $c = 0$, for then $P(x)/Q(x)$ is a scalar multiple of $x^k$ which is impossible for $\deg P < i < k$. Also, if either $P$ or $Q$ has degree less than $i - 1$, then the left-hand side of \eqref{eq:xPQc} has degree less than $2 i - 2$, which is impossible.

  Suppose $\Delta_{k,i+1} = 0$ also. Then there exist rational polynomials $R$, $S$ of degree $i$ for which
  \begin{equation*}
    \frac{R(x)}{S(x)} \equiv x^k \mod h(x)^{i+1},
  \end{equation*}
  and hence
  \begin{equation*}
    0 \equiv \frac{P(x)}{Q(x)} - \frac{R(x)}{S(x)} = \frac{P(x) S(x) - Q(x) R(x)}{Q(x) S(x)} \mod h(x)^i
  \end{equation*}
  The numerator $P(x) S(x) - Q(x) R(x)$ has degree at most $2i - 1$ and is divisible by the degree $2i$ polynomial $h(x)^i$, hence vanishes. Thus the Pad\'e approximants $P/Q$, $R/S$ are actually equal. Since $\gcd\{R,S\}$ is at most linear and in particular coprime to $h$, we obtain that $P/Q = R/S$ is actually a Pad\'e approximant achieving higher accuracy
  \begin{equation*}
    P(x) \equiv x^k Q(x) \mod h(x)^{i+1}.
  \end{equation*}
  Taking logarithmic derivatives, as above, yields $c h(x)^{i-1} \equiv 0 \mod h(x)^i$ which is impossible.
\end{proof}

\begin{rem}
  Instead of $(u_{n+1}/u_n)^k$, we could have considered the continued fraction expansion of $f(u_{n+1}/u_n)$ for any $f(x) \in \QQ(x)$. In this case the determinants $\Delta_{f,i}(a,b)$ of the corresponding linear systems may vanish for multiple consecutive values of $i$, but invariably in strings of odd length $i = i_0, i_0 + 1,\ldots, i_0 + 2r$, corresponding to $f$ having a Pad\'e approximant of (lowest-terms) degree $i_0 - 1$ and accuracy $O(h(x)^{i_0 + r})$, and a monster term $\asymp \lvert\alpha/\bar\alpha\rvert^{(r + 1)n}$.
\end{rem}

\section{Matrix factorization and patch-junction decomposition} \label{sec:matrix_fzn}
We now restrict to the case that all $\Delta_{k,i}(a,b) \neq 0$, continuing to assume that $b = \pm 1$ and $a \geq 2 - b$. If $R \subseteq \RR$ is a subring, we denote by $\V_{k,i}(R)$ the $R$-submodule of pairs $\begin{bsmallmatrix}
  P \\ Q
\end{bsmallmatrix}$ of polynomials of degree at most $i$ with coefficients in $R$ satisfying the Pad\'e relation
\begin{equation}\label{eq:x_pade}
  P(x) \equiv x^k Q(x) \mod h(x)^i.
\end{equation}
Three rings $R$ will be of use to us: $\ZZ$, $\QQ$, and $K = \QQ(\alpha) = \QQ(\sqrt{D})$, where $D = a^2 + 4b$ is the discriminant of $h$. We first construct a $\QQ$-basis of $\V_{k,i}(\QQ)$.
\begin{lem}\label{lem:Phi} For $0 \leq i \leq k$, there is a unique $2 \times 2$ matrix $\boldsymbol{\Phi}_i(X, Y)$ of homogeneous rational polynomials of degree $i$ satisfying the following two conditions:
  \begin{enumerate}[label=\textup{(\roman*)}]
    \item\label{exists:pade} Each column of the dehomogenization $\boldsymbol{\Phi}_i(x, 1)$ belongs to $\V_{k,i}(\QQ)$, in other words,
    \begin{equation} \label{eq:Pade_vec_mat}
      \begin{bmatrix}
        -1 & x^k
      \end{bmatrix} \boldsymbol{\Phi}_i(x, 1) \equiv 0 \mod h^i.
    \end{equation}
    \item\label{exists:ldg} $\boldsymbol{\Phi}_i(1, 0) = \mathbf{I}_2$.
  \end{enumerate}
  Moreover, the columns of $\boldsymbol{\Phi}_i(x, 1)$ form a $\QQ$-basis of $\V_{k,i}(\QQ)$, and $\det \boldsymbol{\Phi}_i(X,Y) = h^{\hom}(X, Y)^i$.
\end{lem}
\begin{proof}
  If the leading coefficients $p, q$ of
  \begin{equation*}
    P(x) = p x^i + \check P(x) \textand Q(x) = q x^i + \check Q(x)
  \end{equation*}
  are fixed, the Pad\'e relation \eqref{eq:x_pade} takes the form
  \begin{equation} \label{eq:pade_inhom}
    x^k \big(q x^i + \check Q(x)\big) - \big(p x^i + \check P(x)\big) = \big(q x^{k-i} + \check E(x)\big) h(x)^i
  \end{equation}
  where $\deg \check E < k - i$. Just as in the proof of Lemma \ref{lem:Pade_Toeplitz}, the coefficients of $x^i, x^{i+1}, \ldots, x^{k-1}$ on the left-hand side of \eqref{eq:pade_inhom} are fixed, which imposes $k - i$ linear conditions on the coefficient vector $\check{\mathbf e}$ of $\check E$, namely
  \begin{equation} \label{eq:euv}
    \mathbf{C}_{k, i} \check{\mathbf e} = p \mathbf{u} + q \mathbf{v}, \qquad \mathbf{u} = -(1, 0, \ldots, 0)^\top, \quad \mathbf{v} = -(c_k, c_{k-1}, \ldots, c_{i+1})^\top.
  \end{equation}
  Since $\Delta_{k,i} = \det \mathbf{C}_{k,i} \neq 0$, \eqref{eq:euv} has a unique rational solution for every $p$ and $q$. Applying this for $(p,q) = (1,0)$ and $(0,1)$, there is a unique such $\boldsymbol{\Phi}_i$ and its dehomogenized columns form a basis for $\V_{k,i}(\QQ)$.

  Since the top row of $\boldsymbol{\Phi}_i(x,1)$ is congruent to $x^k$ times the bottom row modulo $h(x)^i$, we have $h(x)^i \mid \det \boldsymbol{\Phi}_i(x, 1)$, and examining the leading term, we have equality.
\end{proof}

\begin{lem}\label{lem:F}~
  \begin{enumerate}[label=\textup{(\alph*)}]
    \item\label{Phi:reln} For $0 \leq i < k$, we have the relation
    \begin{equation} \label{eq:Phi_reln}
      \boldsymbol{\Phi}_{i+1} = \boldsymbol{\Phi}_i \cdot \big((X - aY) \mathbf{I}_2 + Y \mathbf{T}_i\big)
    \end{equation}
    where $\mathbf{T}_i \in \GL_2(\QQ)$ is a matrix of constants whose characteristic polynomial is $h$.
    \item\label{Phi:fzn} For $n \geq 1$ and $0 \leq i \leq k$, we have a matrix factorization of $\boldsymbol{\Phi}_i$ evaluated at $(u_{n+1}, u_n)$:
    \begin{equation} \label{eq:Phi_fzn}
      \mathbf{U}_i(n) = \boldsymbol{\Phi}_i\left(u_{n+1}, u_n\right) = \mathbf{T}_0^n \mathbf{T}_1^n \cdots \mathbf{T}_{i-1}^n.
    \end{equation}
    In particular,
    \begin{equation} \label{eq:Phi_fzn_special}
      \begin{bmatrix}
        u_{n+1}^k \\ u_n^k
      \end{bmatrix} =
      \mathbf{U}_k(n) \begin{bmatrix}
        1 \\ 0
      \end{bmatrix} \\
      = \mathbf{T}_0^n \mathbf{T}_1^n \cdots \mathbf{T}_{k-1}^n \begin{bmatrix}
        1 \\ 0
      \end{bmatrix}.
    \end{equation}
  \end{enumerate}
\end{lem}
\begin{proof}
  \begin{enumerate}[label=\textup{(\alph*)}]
    \item Observe that each column of
    \begin{equation*}
      \boldsymbol{\Psi} = \boldsymbol{\Phi}_{i+1} - x \boldsymbol{\Phi}_i
    \end{equation*}
    belongs to $\V_{k,i}(\QQ)$ since the $x^{i+1}$ terms cancel. Consequently $\boldsymbol{\Psi} = \boldsymbol{\Phi}_i\mathbf{B}$ for some $\mathbf{B} \in \Mat^{2\times 2}\QQ$, in other words,
    \begin{equation} \label{eq:Phi_reln_temp}
      \boldsymbol{\Phi}_{i+1} = \boldsymbol{\Phi}_i \cdot \big(x \mathbf{I}_2 - \mathbf{B}\big).
    \end{equation}
    Taking $\mathbf{T}_i = a \mathbf{I}_2 - \mathbf{B}$ and homogenizing yields \eqref{eq:Phi_reln}. Taking the determinant of both sides of \eqref{eq:Phi_reln_temp} and canceling $h(x)^{i}$ gives
    \begin{equation*}
      h(x) = \det(x \mathbf{I}_2 - \mathbf{B}_i),
    \end{equation*}
    that is, the characteristic polynomial of $\mathbf{B}_i$ is $h$. Then $\mathbf{T}_i = a\mathbf{I}_2 - \mathbf{B}_i = -b\mathbf{B}_i^{-1}$ is the adjugate of $\mathbf{B}_i$ and thus shares the same characteristic polynomial.

    \item
    By recursively applying \eqref{eq:Phi_reln}, starting from $\boldsymbol{\Phi}_0 = \mathbf{I}_2$, we have
    \begin{equation}
      \boldsymbol{\Phi}_i = \big((X - a Y) \mathbf{I}_2 + Y \mathbf{T}_0\big)\big((X - a Y) \mathbf{I}_2 + Y \mathbf{T}_1\big) \cdots \big((X - a Y) \mathbf{I}_2 + Y \mathbf{T}_{i-1} \big),
    \end{equation}
    so
    \begin{equation}
      \begin{aligned}
        \mathbf{U}_i(n) &= \boldsymbol{\Phi}_i\left(u_{n+1}, u_n\right) \\
        &= \big((u_{n+1} - a u_n) \mathbf{I}_2 + u_n \mathbf{T}_0\big)\big((u_{n+1} - a u_n) \mathbf{I}_2 + u_n \mathbf{T}_1\big) \cdots \big((u_{n+1} - a u_n) \mathbf{I}_2 + u_n \mathbf{T}_{i - 1}\big) \\
        &= \big(b u_{n-1} \mathbf{I}_2 + u_n \mathbf{T}_0\big) \big(b u_{n-1} \mathbf{I}_2 + u_n \mathbf{T}_1\big) \cdots \big(b u_{n-1} \mathbf{I}_2 + u_n \mathbf{T}_{i-1}\big)
      \end{aligned}
    \end{equation}
    It remains to show that
    \begin{equation} \label{eq:T_powers}
      \mathbf{T}_i^n = b u_{n-1}\mathbf{I}_2 + u_n \mathbf{T}_i.
    \end{equation}
    By the Cayley--Hamilton theorem, the powers of $\mathbf{T}_i$ satisfy a linear recurrence with characteristic polynomial $h$, and since the right-hand side of \eqref{eq:T_powers} satisfies the same recurrence, it suffices to check two base cases, say $n = 0$ and $n = 1$. (In the $n = 0$ case we have $u_{-1} = 1/b$, by extending the recurrence backward.)

    To prove \eqref{eq:Phi_fzn_special}, take $i = k$ and note that the left column of $\boldsymbol{\Phi}_k$ is necessarily $\begin{bsmallmatrix}
      X^k \\ Y^k
    \end{bsmallmatrix},$ as $x^k = x^k/1$ is its own Pad\'e approximant of order $k$, the unique one with leading coefficient vector $\begin{bsmallmatrix}
      1 \\ 0
    \end{bsmallmatrix}$.
    \qedhere
  \end{enumerate}
\end{proof}

\subsection{The eigenvectors of \texorpdfstring{$\mathbf{T}_i$}{Tᵢ}}\label{sec:eig}

Each $\mathbf{T}_i$ has characteristic polynomial $h$ and therefore is diagonalizable over $K$. We have the following useful description of the eigenvectors:
\begin{lem}\label{lem:extra_contact}
  Let $\boldsymbol{\xi}_i = \begin{bsmallmatrix}
    \xi_i \\ 1
  \end{bsmallmatrix}$ be the eigenvector of $\mathbf{T}_i$ with eigenvalue $\alpha$. Then
  \begin{equation*}
    \begin{bmatrix}
      P(x) \\ Q(x)
    \end{bmatrix} = \boldsymbol{\Phi}_i(x, 1) \boldsymbol{\xi}_i \in \V_{k,i}(K)
  \end{equation*}
  is the unique solution, up to scaling, of the strengthened Pad\'e condition
  \begin{equation} \label{eq:pade_extra}
    P(x) \equiv x^k Q(x) \mod (x - \alpha)^{i+1} (x - \bar\alpha)^i
  \end{equation}
  with one degree of extra contact at $\alpha$.
\end{lem}
\begin{proof}
  If $\boldsymbol{\xi}_i$ has eigenvalue $\alpha$, then
  \begin{equation*}
    \begin{aligned}
      \V_{k,i+1}(K) &\ni \boldsymbol{\Phi}_{i+1}(x, 1) \boldsymbol{\xi}_i \\
      &= \boldsymbol{\Phi}_i(x, 1)\, \big((x - a) \mathbf{I}_2 + \mathbf{T}_i\big) \boldsymbol{\xi}_i \\
      &= (x - a + \alpha) \boldsymbol{\Phi}_i(x, 1) \boldsymbol{\xi}_i \\
      &= (x - \bar\alpha) \boldsymbol{\Phi}_i(x, 1) \boldsymbol{\xi}_i;
    \end{aligned}
  \end{equation*}
  in other words,
  \begin{equation*}
    (x - \bar\alpha) P(x) \equiv x^k (x - \bar\alpha) Q(x) \mod h(x)^{i+1}.
  \end{equation*}
  Dividing out the common factor $x - \bar\alpha$ yields \eqref{eq:pade_extra}. The solution to \eqref{eq:pade_extra} is unique up to scaling unless \eqref{eq:pade_extra} is satisfied for all $\begin{bsmallmatrix}
    P \\ Q
  \end{bsmallmatrix} \in \V_{k,i}(K)$, in which case the conjugate condition
  \begin{equation*}
    P(x) \equiv x^k Q(x) \mod (x-\alpha)^i (x - \bar\alpha)^{i+1}
  \end{equation*}
  is also satisfied on all of $\V_{k,i}(K)$, and we obtain $\V_{k,i}(K) \subseteq \V_{k,i+1}(K)$, which is impossible.
\end{proof}

Set
\begin{equation*}
  \eta_i(n) = \mathbf{U}_i(n) \xi_i = \boldsymbol{\Phi}_i(u_{n+1}, u_n) \xi_i \in \PP^1(K), \qquad 0 \leq i < k,
\end{equation*}
and $\eta_k(n) = (u_{n+1}/u_n)^k$.
\begin{rem}
  Thanks to the extra contact near $x = \alpha$, we have $\eta_i(n) = (u_{n+1}/u_n)^k + O(\alpha^{-2(i+1)n})$ as $n \to \infty$. The significance of $\eta_i(n)$ is that it is an approximation that ``gets stuck in one patch'': for $a = b = 1$ and $k = 3$, compare \eqref{eq:fib_3} with
  \begin{equation}
    \begin{aligned}
      \eta_0(11) &= [\overline{4}] \\
      \eta_1(11) &= [4,4,4,\mathbf{10},\overline{4}] \\
      \eta_2(11) &= [4,4,4,\mathbf{10},4,4,\mathbf{2},\mathbf{2}, \overline{1}] \\
      (F_{12}/F_{11})^3 = \eta_3(11) &= [4,4,4,\mathbf{10},4,4,\mathbf{2},\mathbf{2},1,1,1,1,1,1,1,1].
    \end{aligned}
  \end{equation}
  We will make our junction and patch words from the preperiodic and periodic parts of the continued fraction expansion of $\eta_i(n)$, respectively.
\end{rem}

\subsection{The integral basis} \label{sec:int_basis}
Within the $\QQ$-vector space $\V_{k,i}(\QQ)$ is the $\ZZ$-lattice $\V_{k,i}(\ZZ)$ of integer-coefficient pairs. Let $\Lambda_{k,i} \subseteq \ZZ^2$ be the set of leading-coefficient pairs of $\V_{k,i}(\ZZ)$, that is, the set of pairs $\begin{bsmallmatrix}
  p \\ q
\end{bsmallmatrix}$ such that $\boldsymbol{\Phi}_i \begin{bsmallmatrix}
  p \\ q
\end{bsmallmatrix}$ has integer coefficients. We note the following:

\begin{lem} \label{lem:Delta_dvd}
  $\Lambda_{k,i} \subseteq \ZZ^2$ is a lattice whose covolume $\widetilde\Delta_{k,i} = [\ZZ^2 : \Lambda_{k,i}]$ divides $\Delta_{k,i}$.
\end{lem}
\begin{proof}
  With notation as in \eqref{eq:euv}, we have $\begin{bsmallmatrix}
    p \\ q
  \end{bsmallmatrix} \in \Lambda_{k,i}$ exactly when $p\mathbf{u} + q\mathbf{v}$ lies in the $\ZZ$-span $\Lambda$ of the columns of $\mathbf{C}_{k,i}$. We therefore have an injection of groups
  \begin{equation*}
    \begin{aligned}
      \ZZ^2 / \Lambda_{k,i} &\longto \ZZ^{k-i} / \Lambda \\
      \begin{bmatrix}
        p \\ q
      \end{bmatrix} & \longmapsto p\mathbf{u} + q\mathbf{v},
    \end{aligned}
  \end{equation*}
  showing that $\widetilde\Delta_{k,i} = [\ZZ^2 : \Lambda_{k,i}]$ is finite and divides $[\ZZ^{k-i} : \Lambda] = |\Delta_{k,i}|$.
\end{proof}

\begin{lem}\label{lem:G_integral}
  If $\widetilde\Delta_{k,i} \mid u_m$ (and, in particular, if $\Delta_{k,i} \mid u_m$), then
  \begin{equation*}
    \mathbf{G}_i(n,m) = \mathbf{U}_i(n + m) \mathbf{U}_i(n)^{-1} \in \GL_2 \ZZ.
  \end{equation*}
\end{lem}
\begin{rem}
  Since $\{u_m\}$ is given by a reversible recurrence, it is purely periodic modulo any $d > 0$, and so there is always a least $m > 0$ for which $d \mid u_m$, commonly called the \emph{rank of apparition} $\rho(d)$. In the introduction, we set $m_i = \rho(\Delta_{k,i})$; we may set $m_i = \rho(\widetilde{\Delta}_{k,i})$ instead, thereby only strengthening the result.
\end{rem}
\begin{proof}[Proof of Lemma \ref{lem:G_integral}]
  Let $\boldsymbol{\Gamma}_i \in \Mat^{2\times 2} \ZZ$ be a matrix whose columns form a $\ZZ$-basis of $\Lambda_{k,i}$. Then $\det \boldsymbol{\Gamma}_i = \pm \widetilde\Delta_{k,i}$. The matrix
  \begin{equation*}
    \widetilde{\boldsymbol{\Phi}}_i = \boldsymbol{\Phi}_i \boldsymbol{\Gamma}_i
  \end{equation*}
  has the property that the columns of its dehomogenization $\widetilde{\boldsymbol{\Phi}}_i(x,1)$ form a $\ZZ$-basis of $\V_{k,i}(\ZZ)$, and
  \begin{equation*}
    \det(\widetilde{\boldsymbol{\Phi}}_i) = \pm \widetilde\Delta_{k,i} \cdot (h^{\hom})^i.
  \end{equation*}
  Then
  \begin{equation*}
    \widetilde{\mathbf{U}}_i(n) = \mathbf{U}_i(n) \boldsymbol{\Gamma}_i = \widetilde{\boldsymbol{\Phi}}_i(u_{n+1}, u_n)
  \end{equation*}
  has integral entries and determinant $\pm \bigl(h^{\hom}(u_{n+1}, u_n)\bigr)^i \widetilde\Delta_{k,i} = \pm \widetilde\Delta_{k,i}$.

  If $\widetilde\Delta_{k,i} \mid u_m$, then $m$ is a period for $\{u_n\}$, up to scaling:
  \begin{equation*}
    u_{n+m} \equiv u_{m+1}u_n \mod \widetilde\Delta_{k,i}.
  \end{equation*}
  Since the entries of $\widetilde{\boldsymbol{\Phi}}_i$ are homogeneous polynomials of degree $i$,
  \begin{equation*}
    \widetilde{\mathbf{U}}_i(n + m) \equiv u_{m+1}^i \widetilde{\mathbf{U}}_i(n) \mod \widetilde\Delta_{k,i}
  \end{equation*}
  for all $n$. Writing
  \begin{equation*}
    \widetilde{\mathbf{U}}_i(n + m) = u_{m+1}^i \widetilde{\mathbf{U}}_i(n) + \widetilde\Delta_{k,i} \mathbf{A}
  \end{equation*}
  for some $\mathbf{A} \in \Mat^{2\times 2} \ZZ$, we get
  \begin{equation*}
    \mathbf{G}_i(n,m) = \widetilde{\mathbf{U}}_i(n + m) \widetilde{\mathbf{U}}_i(n)^{-1} = u_{m+1}^i \mathbf{I}_2 + \widetilde\Delta_{k,i} \mathbf{A} \big(\widetilde{\mathbf{U}}_i(n)\big)^{-1},
  \end{equation*}
  as desired.
\end{proof}

\subsection{Continued fraction factorization}
The connection between continued fraction words and matrices is well known. Here we use the conventions of Mercat \cite{Mercat13}, up to an inconsequential conjugation by $\begin{bsmallmatrix} 0 & 1 \\ 1 & 0
\end{bsmallmatrix}$.
Let $\W$ be the submonoid of $\GL_2\ZZ$ generated by the matrices
\begin{equation}
  [a] = \begin{bmatrix}
    a & 1 \\ 1 & 0
  \end{bmatrix}, \qquad a \in \ZZ_{>0}.
\end{equation}
Then each $\mathbf{W} \in \W$ can be represented by a continued fraction word: $\mathbf{W} = [a_1, a_2, \ldots, a_n] = [a_1][a_2]\cdots [a_n]$. Using the uniqueness of continued fraction expansions, it is easy to see that $\W$ is freely generated \cite[Remarque 2.6]{Mercat13}, that is, each $\mathbf{W} \in \W$ corresponds to a unique continued-fraction word $[a_1, a_2, \ldots, a_n]$, and we freely conflate a word with its corresponding matrix. Thus a word of integers $[a_1, \ldots, a_n]$ can denote one of three objects: a continued fraction, a matrix in $\GL_2 \ZZ$, or a M\"obius transformation in $\PGL_2 \ZZ$ acting on the projective line by $x \mapsto [a_1; a_2,\ldots, a_n, x]$. The concatenation operations on these three types of objects are compatible, and context will determine how a word is being used.

If $\sigma \geq 1$ is a real number, then a continued fraction expansion of $\sigma$ begins with $\mathbf{W}$ if and only if $\mathbf{W}^{-1}(\sigma) \in [1,\infty]$; and if $\mathbf{W}^{-1}(\sigma) = \infty$, then $\mathbf{W}$ is one of the two continued-fraction expansions of $\sigma$. The $\sigma \in \RR$ with purely periodic continued fractions $\sigma = [\overline{\mathbf{W}}]$ are, as is well known, those quadratic irrationals $\sigma > 1$ for which $-1 < \bar\sigma < 0$.

We have the following simple criterion for when a matrix arises from a continued-fraction word:
\begin{prop}
  A matrix $\mathbf{A} \in \GL_2 \ZZ$, $\mathbf{A} \neq \mathbf{I}_2$, belongs to $\W$ if and only if the entries of $\mathbf{A}$ are nonnegative and each row and column of $\mathbf{A}$ is weakly decreasing.
\end{prop}
\begin{proof}
  The proof is essentially by performing the Euclidean algorithm; see Mercat \cite[Proposition 2.3]{Mercat13}.
\end{proof}

We are now ready to construct the patch and junction words and finish the proof of Theorem \ref{thm:main_intro}\ref{main:JPJ}.

\subsubsection{Constructing the patches}

Recall from Section \ref{sec:eig} that
\begin{equation*}
  \eta_i(n) = \mathbf{U}_i(n) \xi_i = \boldsymbol{\Phi}_i(u_{n+1}, u_n) \xi_i.
\end{equation*}
Note that
\begin{equation}\label{eq:G_i}
  \mathbf{G}_i(n, m) \eta_i(n) = \mathbf{G}_i(n, m) \mathbf{U}_i(n) \xi_i = \mathbf{U}_i(n + m) \xi_i = \eta_i(n + m).
\end{equation}
By Lemma \ref{lem:G_integral}, if $m_i \mid m$, then $\mathbf{G}_i(n, m) \in \GL_2\ZZ$ so the periodic part of $\eta_i(n)$ depends only on $n$ modulo $m_i$. On the other hand, since $\xi_i$ is an eigenvector for $\mathbf{T}_i$,
\begin{equation}\label{eq:G_i+1}
  \mathbf{G}_{i+1}(n, m) \eta_i(n) = \mathbf{G}_{i+1}(n, m) \mathbf{U}_i(n) \mathbf{T}_i^n \xi_i = \mathbf{G}_{i+1}(n, m) \mathbf{U}_{i+1}(n) \xi_i = \mathbf{U}_{i+1}(n+m)\xi_i = \eta_i(n+m).
\end{equation}
Consequently, the periodic part $\mathbf{P}_{i,n}$ of the continued fraction for $\eta_i(n)$ also depends only on $n$ modulo $m_{i+1}$, and hence only on $n$ modulo $\gcd(m_i,m_{i+1}) = g_i$. Let $\mathbf{P}_{i,n}$ be the least period of the continued fraction for $\eta_i(n)$, which is unique up to cyclically permuting the terms and only depends on $n$ modulo $g_i$. This will serve as our patch word.

\subsubsection{Constructing the junctions}

For $0 \leq i < k$, let $\sigma_{i,n}$ be the purely periodic real whose continued fraction expansion is $[\overline{\mathbf{P}_{i,n}}]$. Then for each $n$, there is a matrix $\mathbf{H}_i(n) \in \GL_2\ZZ$ (a ``half-junction'') for which
\begin{equation} \label{eq:H_property}
  \mathbf{H}_i(n) \sigma_{i,n} = \eta_i(n),
\end{equation}
unique up to negation and right-multiplication by powers of $\mathbf{P}_{i,n \mod g_i}$. There are two natural families of normalizations of $\mathbf{H}_i(n)$. If $m_i \mid m$, then by \eqref{eq:G_i},
\begin{equation*}
  \eta_{i}(n + m) = \mathbf{G}_i(n, m) \eta_i(n) = \mathbf{G}_i(n, m) \mathbf{H}_i(n) \sigma_{i,n},
\end{equation*}
so we may fix $\mathbf{H}_i^+(n)$ for $n$ in a residue system modulo $m_i$ and extend via
\begin{equation} \label{eq:H-}
  \mathbf{H}_i^+(n + m) = \mathbf{G}_i(n, m) \mathbf{H}_i^+(n), \qquad m_i \mid m.
\end{equation}
Analogously, applying \eqref{eq:G_i+1}, we may fix $\mathbf{H}_i^-(n)$ for $n$ in a residue system modulo $m_{i+1}$ and extend via
\begin{equation*}
  \mathbf{H}_i^-(n + m) = \mathbf{G}_{i+1}(n, m) \mathbf{H}_i^-(n), \qquad m_{i+1} \mid m.
\end{equation*}
For $0 < i < k$, let
\begin{equation*}
  \mathbf{J}_{i, n} = \mathbf{H}_{i-1}^-(n)^{-1} \mathbf{H}_{i}^+(n) \in \GL_2\ZZ,
\end{equation*}
noting that if $m_i \mid m$, then
\begin{equation} \label{eq:junctions_periodic}
  \mathbf{J}_{i, n+m} = \mathbf{H}_{i-1}^-(n)^{-1} \mathbf{G}_i(n, m)^{-1} \mathbf{G}_i(n, m) \mathbf{H}_{i}^+(n) = \mathbf{J}_{i, n},
\end{equation}
so $\mathbf{J}_{i, n}$ depends only on $n$ modulo $m_i$.

It is not obvious that $\mathbf{J}_{i,n} \in \GL_2\ZZ$ is a continued fraction word. However, the following lemma shows that we can resolve the ambiguities in the half-junctions $\mathbf{H}_i^\pm(n)$ to fulfill the needed inequalities.

\begin{lem} Let $\sigma$ and $\tau$ be two reals with purely periodic continued fraction expansions $[\overline{\mathbf{P}}]$, $[\overline{\mathbf{Q}}]$. If $\mathbf{J} \in \GL_2 \ZZ$ and $\mathbf{J}\tau \neq \bar\sigma$, then $\mathbf{P}^m \mathbf{J} \mathbf{Q}^n \in \pm\W$ for all sufficiently large $m$ and $n$.
\end{lem}
\begin{proof}
  Diagonalize
  \begin{equation*}
    \mathbf{P} = \begin{bmatrix}
      w_{11} & w_{12} \\
      w_{21} & w_{22}
    \end{bmatrix}
    =
    \begin{bmatrix}
      \sigma & \bar\sigma \\
      1 & 1
    \end{bmatrix}
    \begin{bmatrix}
      \lambda_\sigma &  \\
      & \lambda_{\bar\sigma}
    \end{bmatrix}
    \begin{bmatrix}
      \sigma & \bar\sigma \\
      1 & 1
    \end{bmatrix}^{-1},
  \end{equation*}
  and note that $\lambda_\sigma = w_{21} \sigma + w_{22} > 1$ and $|\lambda_{\bar\sigma}| = 1/\lambda_\sigma < 1$. Consequently,
  \begin{equation}
    \lim_{m \to \infty} \frac{1}{\lambda_\sigma^m} \mathbf{P}^m =  \begin{bmatrix}
      \sigma & \bar\sigma \\
      1 & 1
    \end{bmatrix}
    \begin{bmatrix}
      1 &   \\
      & 0
    \end{bmatrix}
    \begin{bmatrix}
      \sigma & \bar\sigma \\
      1 & 1
    \end{bmatrix}^{-1}
    = \frac{1}{\sigma - \bar\sigma} \begin{bmatrix}
      \sigma \\ 1
    \end{bmatrix} \begin{bmatrix}
      1 & -\bar\sigma
    \end{bmatrix}.
  \end{equation}
  Treating $\mathbf{Q}$ analogously, we get
  \begin{equation}
    \lim_{m,n\to \infty} \frac{1}{\lambda_\sigma^m \lambda_\tau^n}  \mathbf{P}^m \mathbf{J} \mathbf{Q}^n = \frac{1}{(\sigma - \bar\sigma)(\tau - \bar\tau)} \begin{bmatrix}
      \sigma \\ 1
    \end{bmatrix} \begin{bmatrix}
      1 & -\bar\sigma
    \end{bmatrix} \mathbf{J} \begin{bmatrix}
      \tau \\ 1
    \end{bmatrix} \begin{bmatrix}
      1 & -\bar\tau
    \end{bmatrix}
    = c \begin{bmatrix}
      \sigma \\ 1
    \end{bmatrix} \begin{bmatrix}
      1 & -\bar\tau
    \end{bmatrix},
  \end{equation}
  where
  \begin{equation*}
    c = \frac{1}{(\sigma - \bar\sigma)(\tau - \bar\tau)}\begin{bmatrix}
      1 & -\bar\sigma
    \end{bmatrix} \mathbf{J} \begin{bmatrix}
      \tau \\ 1
    \end{bmatrix}
  \end{equation*}
  is a real number, nonzero because $\mathbf{J}\tau \neq \bar\sigma$. Negating $\mathbf{J}$ if necessary, we may assume that $c > 0$. Then for sufficiently large $m$ and $n$, $\mathbf{P}^m \mathbf{J} \mathbf{Q}^n$ has weakly decreasing rows and columns, because
  \begin{equation*}
    \begin{bmatrix}
      \sigma \\ 1
    \end{bmatrix} \begin{bmatrix}
      1 & -\bar\tau
    \end{bmatrix}
  \end{equation*}
  does, and thus lies in $\W$.
\end{proof}
To apply this lemma with
\begin{equation*}
  \mathbf{J} = \mathbf{J}_{i, n} = \mathbf{H}_{i-1}^-(n)^{-1} \mathbf{H}_{i}^+(n), \qquad \tau = \sigma_{i, n}, \qquad \sigma = \sigma_{i-1, n},
\end{equation*}
we must show that $\mathbf{J}_{i,n}\sigma_{i,n} \neq \bar\sigma_{i-1,n}$. Suppose the contrary. Then
\begin{equation*}
  \begin{aligned}
    \xi_i &= \mathbf{U}_i(n)^{-1} \eta_i(n) \\
    &= \mathbf{U}_i(n)^{-1} \mathbf{H}_i^+(n) \sigma_{i, n} \\
    &= \mathbf{U}_i(n)^{-1} \mathbf{H}_{i-1}^-(n) \mathbf{J}_{i,n} \sigma_{i,n} \\
    &= \mathbf{U}_i(n)^{-1} \mathbf{H}_{i-1}^-(n) \bar\sigma_{i-1,n} \\
    &= \mathbf{U}_i(n)^{-1} \bar\eta_{i-1}(n) \\
    &= \mathbf{T}_{i-1}^{-n} \mathbf{U}_{i-1}(n)^{-1} \bar\eta_{i-1}(n) \\
    &= \mathbf{T}_{i-1}^{-n} \bar\xi_{i-1} \\
    &= \bar\xi_{i-1}.
  \end{aligned}
\end{equation*}
In other words, $\mathbf{T}_i$ and $\mathbf{T}_{i-1}$ have the same eigenvectors $\boldsymbol{\xi}_i, \bar{\boldsymbol{\xi}}_i$ but with interchanged eigenvalues, and $\mathbf{T}_{i-1} = a\mathbf{I}_2 - \mathbf{T}_i$. We obtain
\begin{equation*}
  \begin{aligned}
    \boldsymbol{\Phi}_{i+1}(X,Y)
    &=\boldsymbol{\Phi}_{i-1}(X,Y)
    \cdot \bigl((X - a Y)\mathbf{I}_2 + Y (a \mathbf{I}_2 - \mathbf{T}_i)\bigr)
    \cdot \bigl((X - a Y)\mathbf{I}_2 + Y (\mathbf{T}_i)\bigr)\\
    &=h^{\hom}(X,Y)\boldsymbol{\Phi}_{i-1}(X,Y).
  \end{aligned}
\end{equation*}
Since $\boldsymbol{\Phi}_{i+1}$ satisfies the Pad\'e relation \eqref{eq:Pade_vec_mat} to order $i + 1$, the columns of $\boldsymbol{\Phi}_{i-1}$ are Pad\'e approximants of degree at most $i - 1$ with contact order at least $i$, so $\Delta_{k,i} = 0$, a contradiction.

Thus, if the half-junctions $\mathbf{H}_i^-(n), \mathbf{H}_i^+(n)$ are right-multiplied by sufficiently high negative and positive powers, respectively, of the appropriate patch word $\mathbf{P}_{i,n}$, and possibly negated, then the junctions $\mathbf{J}_{i,n}$ eventually land in $\W$ (since $\mathbf{J}_{i, n+m_i} = \mathbf{J}_{i,n}$, only finitely many values of $n$ need be considered).

The beginning and end of the continued fraction deserve special attention. For the first junction, Lemma \ref{lem:extra_contact} gives $\eta_0(n) = \alpha^k$, and we set $\mathbf{J}_{0,n} = \mathbf{H}_0^+(n) = \mathbf{H}_0^+(0)$ to be the preperiodic part of the continued fraction of $\alpha^k$, the periodic part being $\mathbf{P}_{0,n} = \mathbf{P}_{0,0}$ (which, as the reader can verify without much effort, consist of at most one and two terms respectively).

As to the last patch, the extra-contact approximation can be derived by solving
\begin{equation}
  x^k Q(x) - P(x) = (x - \alpha)^k (x - \bar\alpha)^{k-1} = h(x)^{k-1} (x - \alpha).
\end{equation}
The leading coefficients of $Q$ and $P$ are respectively $1$ and $c \alpha - d$, where $c = c_{k-1}$ and $d = c_k$ are coefficients of $h(x)^{k-1}$. Note that $c = \Delta_{k,k-1} \neq 0$. We obtain $\eta_{k-1}(0) = \xi_{k-1} = (P/Q)(\infty) = c \alpha - d$, which, like any algebraic integer of degree $2$, is ``almost purely periodic'': we may take
\begin{equation*}
  \sigma_{k-1,0} = |c| \alpha - \ceil{\strut|c| \bar\alpha},
\end{equation*}
the preimage of $\eta_{k-1}(0)$ under the unimodular matrix
\begin{equation*}
  \mathbf{H}_{k-1}^-(0) = \begin{bmatrix}
    \sgn c & (\sgn c) \ceil{\strut|c| \bar\alpha} - d \\
    0 & 1
  \end{bmatrix} \in \GL_2 \ZZ
\end{equation*}
which, usefully, fixes $\infty$. Since $m_k = 1$, this choice determines $\mathbf{H}_{k-1}^-(n) = \mathbf{G}_k(0,n) \mathbf{H}_{k-1}^-(0)$ for all $n$.

\begin{rem}
  If we were writing this proof for formalized verification, we would have fixed $\mathbf{P}_{k-1,n}$ and its associated data first, then $\mathbf{P}_{0,n}$ which must be taken to be the same $\mathbf{P}_{k-1,n}$ in the edge case $k = 1$, and lastly the generic $\mathbf{P}_{i,n}$. But this paper is not a summary of a formalization, and the case $k = 1$ is hardly the most motivating.
\end{rem}

\subsubsection{Putting it together}
\begin{proof}[Proof of Theorem \ref{thm:main_intro}\ref{main:JPJ}]
  The final step is to open up the matrix factorization \eqref{eq:Phi_fzn_special} to reveal the continued fraction structure. Let $n_i$ be the least nonnegative residue of $n$ modulo $m_i$. Each $\mathbf{T}_i^n$ factor may be written as follows:
  \begin{equation} \label{eq:T_6_factors}
    \begin{aligned}
      \mathbf{T}_i^n &= \mathbf{U}_i(n)^{-1} \mathbf{U}_{i+1}(n) \nonumber \\
      &= \mathbf{U}_i(n_i)^{-1} \mathbf{G}_i(n_i, n - n_i)^{-1} \mathbf{G}_{i+1}(n_{i+1}, n - n_{i+1}) \mathbf{U}_{i+1}(n_{i+1}) \nonumber \\
      &= \mathbf{U}_i(n_i)^{-1} \mathbf{H}_i^+(n_i) \mathbf{H}_i^+(n)^{-1} \mathbf{H}_i^-(n) \mathbf{H}_i^-(n_{i+1})^{-1} \mathbf{U}_{i+1}(n_{i+1}) \\
      &= \pm \mathbf{U}_i(n_i)^{-1} \mathbf{H}_i^+(n_i) \mathbf{P}_{i,n}^{e_i(n)} \mathbf{H}_i^-(n_{i+1})^{-1} \mathbf{U}_{i+1}(n_{i+1})
    \end{aligned}
  \end{equation}
  where $e_i(n)$ is an integer, since the only elements of $\GL_2\ZZ$ fixing $\sigma_{i, n}$ are $\pm \mathbf{P}_{i,n}^e$ for $e \in \ZZ$. If $n$ grows within a fixed congruence class modulo $\lcm\{m_i, m_{i+1}\}$, then only the middle factor $\mathbf{P}_{i,n}^{e_i(n)}$ can change, and the growth of $e_i(n)$ is controlled by the eigenvalue:
  \begin{equation}
    \begin{aligned}
      \alpha^n \begin{bmatrix}\xi_i \\ 1\end{bmatrix}
      &= \mathbf{T}_i^n \begin{bmatrix}\xi_i \\ 1\end{bmatrix} \\
      &= c \mathbf{U}_i(n_i)^{-1} \mathbf{H}_i^+(n_i) \mathbf{P}_{i,n}^{e_i(n)} \cdot \begin{bmatrix}
        \sigma_{i, n_{i+1}} \\ 1
      \end{bmatrix} \\
      &= c \lambda^{e_i(n)} \mathbf{U}_i(n_i)^{-1} \mathbf{H}_i^+(n_i)  \begin{bmatrix}
        \sigma_{i, n_{i+1}} \\ 1
      \end{bmatrix} \\
      &= c' \lambda^{e_i(n)} \begin{bmatrix}\xi_i \\ 1\end{bmatrix}
    \end{aligned}
  \end{equation}
  where $c,c'$ are constants and $\lambda > 1$ is the dominant eigenvalue of the patch $\mathbf{P}_{i,n}$. Consequently we have the relation
  \begin{equation*}
    \alpha^n = c' \lambda^{e_i(n)},
  \end{equation*}
  which shows that $e_i(n)$ is linear in $n$ in this congruence class, and also nonnegative for $n$ sufficiently large.

  To finish, we multiply together all the $\mathbf{T}_i$:
  \begin{equation}
    \begin{aligned}
      \left(\frac{u_{n+1}}{u_n}\right)^k &=
      \left(\prod_{i=0}^{k-1} \mathbf{T}_i^n\right) (\infty) \\
      &= \left(\prod_{i = 0}^{k-1} \mathbf{U}_i(n_i)^{-1} \mathbf{H}_i^+(n_i) \mathbf{P}_{i,n_i}^{e_i(n)} \mathbf{H}_i^-(n_{i+1})^{-1} \mathbf{U}_{i+1}(n_{i+1}) \right)(\infty) \\
      &= \mathbf{U}_0(0)^{-1} \mathbf{H}_0^{+}(0) \mathbf{P}_{0,0}^{e_0(n)} \left(\prod_{i=1}^{k-1} \mathbf{H}_{i-1}^-(n_i)^{-1} \mathbf{U}_{i}(n_i) \mathbf{U}_{i}(n_i)^{-1} \mathbf{H}_{i}^+(n_i) \mathbf{P}_{i,n_i}^{e_i(n)}\right) \mathbf{H}_{k-1}^-(n_k)^{-1} \mathbf{U}_k(n_k)(\infty) \\
      &= \mathbf{H}_0^+(0) \mathbf{P}_{0,0}^{e_0(n)} \left(\prod_{i=1}^{k-1} \mathbf{H}_{i-1}^-(n_i)^{-1}  \mathbf{H}_i^+(n_i) \mathbf{P}_{i,n_i}^{e_i(n)}\right) \mathbf{H}_{k-1}^-(0)^{-1} (\infty) \\
      &= \mathbf{J}_{0,0} \mathbf{P}_{0,0}^{e_0(n)} \left(\prod_{i=1}^{k-1} \mathbf{J}_{i,n} \mathbf{P}_{i,n}^{e_i(n)} \right) (\infty),
    \end{aligned}
  \end{equation}
  establishing the desired continued fraction decomposition.
\end{proof}

\section{Generalizations and epilogue} \label{sec:epilogue}

The restriction $b = \pm 1$ in this proof was mainly to decrease the amount of bookkeeping in an already rather technical proof. With the same techniques and added bookkeeping, one can show that the same sort of decomposition is possible if $b \mid a^2$.
\begin{exer}\label{thm:main_gen}
  Let $a$ and $b$ be integers with $a > 0$, $D = a^2 + 4 b > 0$ a non-square, and $b \mid a^2$, so that the ratio $\rho = \alpha/\bar\alpha$ of the roots of $h(x) = x^2 - a x - b$ is a real irrational algebraic unit. Consider the corresponding recurrence
  \begin{equation*}
    u_0 = 0, \quad u_1 = 1, \quad u_{n+2} = a u_{n+1} + b u_n.
  \end{equation*}
  \begin{enumerate}[label=\textup{(\alph*)}]
    \item If $\Delta_{k,i}(a,b) = 0$ for some $i$, then for all sufficiently large $n$, the power $(u_{n+1}/u_n)^k$ has a continued fraction with a ``monster term\kern-1pt'' $M$ bounded below and above by constant multiples of $|\rho|^n$.
    \item If, on the other hand, $\Delta_{k,i}(a,b) \neq 0$ for all $i$, then for $1 \leq i \leq k$, let $m_i$ be the least positive even integer $n$ for which
    \begin{equation*}
      |b|^i \widetilde{\Delta}_{k,i} \mathop{\Big|} \frac{u_n}{a b^{n/2 - 1}},
    \end{equation*}
    which exists because the terms $v_n = u_{2n}/(a b^{n - 1})$ satisfy the reversible recurrence
    \begin{equation*}
      v_{n+2} = \left(\frac{a^2}{b} + 2\right) v_{n+1} - v_n.
    \end{equation*}
    Set $m_0 = 1$ and $g_i = \gcd\{m_i,m_{i+1}\}$.

    Then for $n$ sufficiently large, $(u_{n+1}/u_n)^k$ has a continued fraction with a patches-and-junctions decomposition
    \begin{equation*}
      \mathbf{J}_{0,0} \mathbf{P}_{0,0}^{e_0(n)} \mathbf{J}_{1, n \bmod m_1} \mathbf{P}_{1, n \bmod g_1}^{e_1(n)} \mathbf{J}_{2, n \bmod m_2} \cdots \mathbf{J}_{k-1, n \bmod m_{k-1}} \mathbf{P}_{k-1, n \bmod g_{k-1}}^{e_{k-1}(n)} \mathbf{J}_{k, n \bmod m_k}
    \end{equation*}
    where the continued-fraction words $\mathbf{P}_{i,n}$, $\mathbf{J}_{i,n}$ depend only on $n$ modulo $g_i$ and $m_i$ respectively, and the number $e_i(n)$ of repetitions within each patch is a linear function of $n$ when $n$ is in a fixed congruence class modulo $\lcm\{m_i, m_{i+1}\}$. In particular, the terms of the continued fraction remain bounded as $n$ grows.
  \end{enumerate}
\end{exer}

\begin{rem}
  In contrast to Theorem \ref{thm:main_intro}, there may be a nonempty final junction $\mathbf{J}_{k,n}$, and its governing modulus $m_k$ may be greater than $1$.
\end{rem}
\begin{rem}
  The condition $b \mid a^2$ makes the root ratio $\rho = \alpha/\bar\alpha$ an algebraic unit, and thus the powers $\rho^m$ can be eigenvalues for continued fraction patches. If $b \nmid a^2$, then no $\alpha^n\bar\alpha^m$ is an algebraic unit, and the initial patch $\mathbf{P}_{0,n}$ grows at an irrational rate, with no periodically recurring junction following, as in \eqref{eq:chaos}.
\end{rem}
\begin{rem}
  Experimentation seems to show that if the initial terms $u_0, u_1$ are altered, the monster criterion $\Delta_{k,i}(a, b)$ is unchanged but the patch periods and moduli $m_i, g_i$ may increase. We leave this for further investigation.
\end{rem}

\subsection{The hunt for additional monsters}

As $k$ increases, the degrees of the Pad\'e--Toeplitz determinants
$\delta_{k,i}(t)$ increase, and the possibility of an integer root $t$ heuristically becomes rarer. Computations have revealed three more integer solutions, in addition
to the one noted in the introduction \eqref{eq:fib_5_table}, which are listed below
in Table \ref{tab:monster}.

\begin{table}[ht]
  \begin{equation*}
    \begin{array}{c|c|c|c}
      k&i&\delta_{k,i}(t)&t=a^2/b\\
      \hline
      3&2&t-2&2\\
      4&3&t-6&6\\
      5&2&(t-1)(t^2+3t+6)&1\\
      6&3&(t-3)(t^3+3t^2+18t+18)&3
    \end{array}
  \end{equation*}
  \caption{The four known cases in which $(u_{n+1}/u_n)^k$ displays a monster term}
  \label{tab:monster}
\end{table}

Each entry yields a one-parameter family of $(a,b)$. For instance, for the last entry in the table, we may take $a=b=3$ and get a monster of size $M\asymp |\rho|^n = \big((5+\sqrt{21})/2\big)^n$:
\begin{equation*}
  \resizebox{\textwidth}{!}{$\displaystyle\begin{aligned}
    (u_{13}/u_{12})^6 &= [2969,1,3,13,7,1,8,1,1,4,2,1,6,10,1,\mathbf{320340},21,7,2,1,1,23,4,15,1,1,16,1,2] \\
    (u_{14}/u_{13})^6 &= [2969,1,3,13,1,1,1,9,1,3,1,1,8,7,1,6,1,4,\mathbf{1534845},1,1,14,6,2,6,1,1,2,1,8,1,1,1,2,1,56,1,5,7,1,2] \\
    (u_{15}/u_{14})^6 &= [2969,1,3,13,1,1,5,2,3,2,3,1,2,1,8,3,1,30,1,\mathbf{7353886},1,2,2,1,2,4,7,2,6,1,1,2,1,8,1,3,1,5,1,1,16,1,2] \\
    (u_{16}/u_{15})^6 &= [2969,1,3,13,1,1,3,3,3,2,1,3,1,5,1,1,1,3,3,1,1,1,16,1,\mathbf{35234591},1,2,11,1,5,2,59,1,31,16,519,1,1,7,3].
  \end{aligned}$}
\end{equation*}
Further computation up to $k = 23$ has revealed that no further $\delta_{k,i}(t)$ has an integer root, apart from extraneous roots at $t = -1,-2,-3$ which correspond to recurrences with non-real $\alpha$. We accordingly conjecture that the list is complete:
\begin{conj}
  The only $(a,b,k)$ for which the continued fraction of $(u_{n+1}/u_n)^k$ has a monster term are the four cases shown in Table \ref{tab:monster}.
\end{conj}

\subsection{Invisible junctions}
There is another exceptional phenomenon worth noting, not nearly as visible as the monster term but even rarer. Looking carefully at the continued fraction expansions \eqref{eq:fib_4} of Fibonacci fourth powers, observe that although Theorem \ref{thm:main_intro} guarantees a decomposition into four patches, only three are visible, with periods
\begin{equation}
  \mathbf{P}_{0,0} = [1, 5], \qquad \mathbf{P}_{1,0} = [3,1,2,3,2,1,3,33], \qquad \mathbf{P}_{2,0} = [11].
\end{equation}
Actually, $\mathbf{P}_{3,0} = [11]$ as well, and the junction $\mathbf{J}_{3,0}$ separating them is empty or, according to the normalization, consists of another finite string of $11$'s. Tracing through the proof, we find that an invisible junction $\mathbf{J}_{i,n}$ occurs exactly when $\xi_i = \xi_{i-1}$ corresponds to a Pad\'e approximant $P/Q \in K(x)$ with \emph{double} extra contact at one of the two basepoints, instead of just one:
\begin{equation*}
  P(x) \equiv x^k Q(x) \mod (x - \alpha)^{i+1} (x - \bar\alpha)^{i-1}, \qquad \deg P, \deg Q \leq i-1
\end{equation*}
This leads to the following analogue of Lemma \ref{lem:Pade_Toeplitz}.

\begin{exer}
  Let $\beta_j \in K$ be the coefficients of
  \begin{equation*}
    (x - \alpha)^{i+1} (x - \bar\alpha)^{i-1} = h(x)^{i-1} (x - \alpha)^2 = \beta_0 x^{2i} + \beta_1 x^{2i - 1} + \cdots + \beta_{2i}.
  \end{equation*}
  If
  \begin{equation}
    \Upsilon_{k,i}(\alpha) = \begin{vmatrix}
      \beta_i & \beta_{i+1} & \cdots & \beta_{k-1} \\
      \beta_{i-1} & \beta_i & \cdots & \beta_{k-2} \\
      \vdots & \vdots & \ddots & \vdots \\
      \beta_{2i-k+1} & \beta_{2i-k+2} & \cdots & \beta_i
    \end{vmatrix} = 0,
  \end{equation}
  then, upon renormalizing the patch-junction decomposition, $\mathbf{J}_{i,n} = \emptyset$ is the empty word and $\mathbf{P}_{i-1,n} = \mathbf{P}_{i,n}$ for all $n$. Conversely, if $\mathbf{J}_{i,n} = \emptyset$ and $\mathbf{P}_{i-1,n} = \mathbf{P}_{i,n}$ for some $n$, then $\Upsilon_{k,i}(\alpha) = 0$.
\end{exer}

Checking $\Upsilon_{k,i}$ for $k \leq 20$, we found no additional invisible junctions, prompting the conjecture:
\begin{conj}
  The only solution to $\Upsilon_{k,i}(\alpha) = 0$ with $a$, $b$ nonzero integers and $D = a^2 + 4 b > 0$ occurs at $k = 4$, $i = 3$, $a^2 = b$, corresponding to the unique invisible junction in $a^4 (F_{n+1}/F_n)^4$.
\end{conj}

We end with a deliberately vague question.

\begin{qn}
  Why are all junctions so short? In \eqref{eq:fib_4}, for instance, each junction is at most three terms (shown in bold), when allowing patches to begin and end anywhere in the cycle. Experimentally, as $k$ increases, the patch lengths skyrocket while the junctions commonly have only one or two terms and can often be empty (but not invisible, the two neighboring patches being distinguishable).
\end{qn}

\newcommand{\etalchar}[1]{$^{#1}$}

\end{document}